\documentclass{amsart}

\usepackage{amsmath,amsfonts,amsthm,amssymb,graphics,graphicx, verbatim, easybmat}

\newtheorem{thm}{Theorem}[section]
\newtheorem{defin}[thm]{Definition}

\newtheorem{lem}[thm]{Lemma}

\newtheorem{prop}[thm]{Proposition}
\theoremstyle{remark} \newtheorem{remark}[thm]{Remark}

\newcommand{\norm}[1]{\left|\!\left|{#1}\right|\!\right|}
\newcommand{\R}{\ensuremath{\mathbb{R}}}
\newcommand{\RR}{\ensuremath{\mathbb{R}}}

\newcommand\bdf{{r}}
\newcommand\xbar{\bar{x}}
\newcommand\wbar{\bar{w}}
\newcommand\vbar{\bar{v}}
\newcommand\xibar{\bar{\xi}}
\newcommand\canon{\mathcal{C}}
\newcommand\Id{\operatorname{Id}}

\title{Semiclassical $L^p$ estimates of quasimodes on curved hypersurfaces}

\author{Andrew Hassell}

\email{Andrew.Hassell@anu.edu.au}

\author{Melissa Tacy}

\email{Melissa.Tacy@anu.edu.au}

\address{Department of Mathematics, Mathematical Sciences Institute, Australian National University, Canberra  0200 ACT AUSTRALIA}


\thanks{This research was supported in part by Australian Research Council Discovery Grant DP0771826, and an Australian Postgraduate Award}

\begin{document}

\begin{abstract}
Let $M$ be a compact manifold of dimension $n$, $P = P(h)$ a semiclassical pseudodifferential operator on $M$, and  $u = u(h)$  an $L^2$ normalised family of functions such that $P u$ is $O(h)$ in $L^2(M)$ as $h \downarrow 0$.  Let $H \subset M$ be a compact  submanifold of $M$. In a previous article, the second-named author proved estimates on the $L^p$ norms, $p \geq 2$, of $u$ restricted to $H$, under the assumption that the $u$ are semiclassically localised and  under some natural structural assumptions about the principal symbol of $P$. These estimates are of the form $C h^{-\delta(n, k, p)}$ where $k = \dim H$ (except for a logarithmic divergence in the case $k = n-2, \, p=2$). When $H$ is a hypersurface, i.e. $k=n-1$, we have $\delta(n, n-1, \, 2) = 1/4$, which is sharp when $M$ is the round $n$-sphere and $H$ is an equator. 

In this article, we assume that $H$ is a hypersurface, and make the additional geometric assumption that $H$ is \emph{curved} (in the sense of Definition~\ref{curvedtoflow} below) with respect to the bicharacteristic flow of $P$. Under this assumption we improve the estimate from $\delta = 1/4$ to $1/6$, generalising work of Burq-G\'erard-Tzvetkov and Hu for Laplace eigenfunctions. To do this we apply the Melrose-Taylor theorem, as adapted by Pan and Sogge, for Fourier integral operators with folding canonical relations. 
\end{abstract}
\maketitle

\section{Introduction}
Let $M$ be a compact manifold of dimension $n$ and $P = P(h)$ a semiclassical pseudodifferential operator on $M$ parametrised by the positive number $h \in (0, h_0]$. Suppose that  $u = u(h)$ is an $O(h)$ quasimode, i.e.  an $L^2$-normalised family of functions, defined for some subset of $(0, h_0]$ accumulating at $0$, such that $P(h) u(h)$ is $O(h)$ in $L^2(M)$.  We assume $P$ has real principal symbol $p(x, \xi)$ and that its full symbol is smooth in $h$. We also put technical assumptions on $p(x,\xi)$ (see Definition \ref{admissible} and \ref{curvedtoflow}) and assume $u$ is localised (see Definition \ref{localised}). One important special case is when $P(h) = h^2 \Delta - 1$ where $\Delta$ is the Laplacian with respect to a Riemannian metric on $M$. Then $u(h)$ is an approximate eigenfunction with eigenvalue $h^{-2}$:
$$
(\Delta - h^{-2}) u(h) = O(h^{-1}) \text{ in } L^2(M).
$$
Other cases of interest are discussed in \cite{koch07}, where this framework was introduced. 

The aim of this paper is to bound the extent to which $u$ can concentrate as $h \to 0$ by estimating the $L^{p}$ norm of $u$ restricted to hypersurfaces, in a manner that is sharp (up to a constant independent of $h$) as $h \to 0$. In particular, we wish to relate the degree of concentration to the geometry of the hypersurface relative to the bicharacteristic flow of $P(h)$.

There are a number of ways to study concentration of eigenfunctions. One can for example study semiclassical measures as in  G\'{e}rard-Leichtnam \cite{gerard}, Zel\-ditch \cite{zelditch}, Zelditch-Zworski \cite{zelditch2}, Anantharaman \cite{anantharaman}, Anantharaman-Koch-Non\-nen\-macher \cite{anantharaman07a}, Anantharaman-Nonnenmacher \cite{anantharaman07}. The aim of these studies is generally to prove non-concentration theorems under geometric conditions on the geodesic flow (such as Anosov flow).

In 1988 Sogge \cite{sogge88} produced sharp $L^{p}$ estimates for spectral clusters (and therefore eigenfunctions) of elliptic operators, comparing the size of the $L^{p}$ norm over the full manifold to the $L^{2}$ norm in terms of powers of the eigenvalue $\lambda$. In 2004 Reznikov \cite{reznikov} proved bounds for restrictions of Laplacian eigenfunctions to curves where the underlying manifold is a hyperbolic surface and in 2007 Burq, G\'{e}rard and Tzvetkov \cite{burq07} proved estimates for general submanifolds and Laplacian eigenfunctions. Their estimates are sharp for sub-sequences of spherical harmonics. For high $p$ these estimates are optimised by eigenfunctions concentrating at a point. For low $p$ the optimising examples are eigenfunctions concentrating in a small tube around a stable periodic geodesic. Burq, G\'{e}rard and Tzvetkov \cite{burq07} were also able to obtain better estimates for small $p$ in dimension two when the submanifold is a curve with positive geodesic curvature. Hu \cite{hu09} extended this to hypersurfaces in $n$ dimensions where the hypersurface has positive curvature. In the special case of a flat two or three dimensional torus Bourgain and Rudnick obtain an improved nonconcentration result for curved hypersurfaces \cite{bourgain09}.
 
In 2009 Tacy \cite{tacy09} extended Burq, G\'{e}rard and Tzvetkov's results on Laplacian eigenfunctions to quasimodes of semiclassical operators. This extension uses the semiclassical framework set up in Burq-G\'{e}rard-Tzvetkov \cite{burq2} and Koch-Tatatru-Zworski \cite{koch07}. The main result of \cite{tacy09} is the following, where we refer to Definitions~\ref{localised} and \ref{admissible} for the precise definitions of localisation and admissibility. 

\begin{thm}\label{tacytheorem}
Let $(M,g)$ be a smooth manifold without boundary and let $H$ be a smooth embedded hypersurface. Let $u(h)$ be a family of $L^{2}$ normalised functions that satisfy $Pu=O_{L^{2}}(h)$ for $P$ a semiclassical operator with symbol $p(x,\xi)$. Assume further that $u$ satisfies the localisation property and that the symbol $p(x,\xi)$ is admissible. Then
$$\norm{u}_{L^{p}(H)}\lesssim{}h^{-\delta(n,p)},$$
\begin{equation}\delta(n,p)=
\begin{cases}
\frac{n-1}{2}-\frac{n-1}{p},&\frac{2n}{n-1}\leq{}p\leq\infty\\
\frac{n-1}{4}-\frac{n-2}{2p},&2\leq{}p\leq\frac{2n}{n-1}
\end{cases}.
\label{delta(n,p)}\end{equation}
\end{thm}

\begin{remark} We have only given the results of \cite{tacy09} pertaining to hypersurfaces. Higher codimension submanifolds were
also treated there. 
\end{remark}

\begin{figure}
 \centering
 \includegraphics[width=9cm,height=7cm]{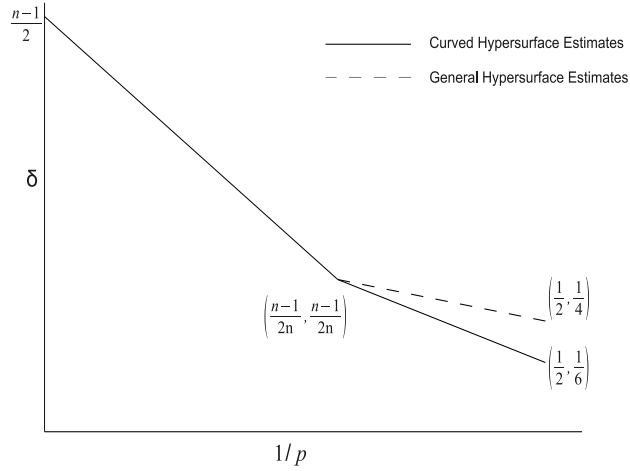}
 \label{fig:hypersurfacebetter}
\caption{$\delta(p)$ plotted against $1/p$ for a general hypersurface  and for a hypersurface curved with respect to the flow.}
\end{figure}

This paper extends the  estimates of Burq-G\'erard-Tzvetkov and Hu for curved hypersurfaces to the semiclassical regime, framing the geometric conditions in terms of the classical (bicharacteristic) flow. 
To motivate the condition of curvature, recall that the classical flow defined by
\begin{equation}\begin{cases}
\dot{x}=\partial_{\xi}p(x,\xi)\\
\dot{\xi}=-\partial_{x}p(x,\xi)\end{cases}\label{classicalflow}\end{equation}
 describes the movement in phase space of a classical particle with classical Hamiltonian $p(x,\xi)$. For the model case of the Laplacian the flow defined by (\ref{classicalflow}) is the geodesic flow. In the semiclassical regime we wish to find estimates that link the properties of this classical flow to concentrations of quasimodes. Intuitively we can think of highly localised packets moving on trajectories defined by the flow. The more time a packet spends near a hypersurface the move concentration we would expect to see there. In \cite{koch07} and \cite{tacy09} it is shown that for a hypersurface $H'$ with boundary defining function\footnote{We say that the real function $r$ is a boundary defining function for $H'$ if $H' = \{ r = 0 \}$ and if $r$ vanishes simply at $H'$, i.e.  $dr \neq 0$ at $H'$.} $\bdf$, if at some point $(x_0, \xi_0)$ we have $\dot \bdf \neq 0$, where the dot indicates derivative with respect to bicharacteristic flow, and if $u$ is a quasimode sufficiently localized near $(x_0, \xi_0)$, then $u$ does not concentrate at $H'$. That is, if $\chi\in{}C_{0}^{\infty}(\R^{n}\times{}\R^{n})$ is a cut off function with small enough support around $(x_{0},\xi_{0})$ then
\begin{equation}
\norm{\chi(x,hD)u}_{L^{2}(H')}\lesssim{}\norm{u}_{L^{2}(M)}.
\label{nonconcentration}\end{equation}
 However, in the general case, a bicharacteristic may stay inside $H$, allowing considerable concentration of an associated wave packet on $H$. As shown in \cite{koch07} and \cite{tacy09}, concentration (as measured by $L^2$ norm) could be as bad as $\sim h^{-1/2}$ assuming just the localisation condition and assumption (A1) below, while additionally assuming (A2) introduces dispersion effects which reduces the concentration to $\sim h^{-1/4}$. To improve on this,
we need to rule out bicharacteristics that stay inside $H$. A natural assumption to make is that the projections of bicharacteristics are  only simply tangent to $H$. In local coordinates this is the same as saying the whenever a bicharacteristic is tangent to $H$, i.e.  $\dot \bdf(x_0, \xi_0)$ vanishes, $x_0 \in H$, then the normal acceleration $\ddot{\bdf}(x_{0},\xi_{0})$ is nonzero. We phrase this by saying that $H$ is \emph{curved} with respect to the bicharacteristic flow. 
 
 Under this additional assumption, which we label (A3) below, we show that the concentration is at most $\sim h^{-1/6}$:

\begin{thm}\label{glancingflowtheorem}
Let $M$, $H$, $P(h)$ and $u(h)$ be as in Theorem~\ref{tacytheorem}. 
If $H$ is  curved with respect to the flow given by $p(x,\xi)$, i.e. satisfies assumption (A3) below, then the estimate \eqref{delta(n,p)} for $p=2$ can be improved from $\delta = 1/4$ to $\tilde \delta = 1/6$. By interpolation with the result for $p = 2n/(n-1)$, we obtain 
\begin{equation}\begin{gathered} \norm{u}_{L^{p}(H)}\lesssim{}h^{-\tilde\delta(n,p)}, \quad 2\leq{}p\leq\frac{2n}{n-1}, \\
\tilde\delta(n,p)=
\frac{n-1}{3}-\frac{2n-3}{3p}, 
\end{gathered}\label{mainestimate}\end{equation}
under assumption (A3). 
\end{thm}

\begin{remark}
For $p\geq{}2n/(n-1)$ there is no improvement in the curved case. 
In this case the $\| \cdot \|_{L^p(H)}$ norm is maximised by functions that concentrate at points so we would not expect the geometry of the hypersurface to affect such estimates. \end{remark}


\section{Semiclassical Analysis}\label{semiclassical}
We work with semiclassical pseudodifferential operators (for a full introduction see \cite{burq2}, \cite{evans} or \cite{koch07}). Such operators are defined by their symbol $p(x,\xi,h)$ and a quantisation procedure
$$P(h) u(h) = p(x,hD, h)u(h)=\frac{1}{(2\pi{}h)^{n}}\int{}e^{\frac{i}{h}<x-y,\xi>}p(x,\xi,h)u(y,h)d\xi{}dy$$
where $h$ is a small parameter. 
Because we are just about to assume that $u$ is localised (Definition \ref{localised}), it is harmless to assume that $p$ is a $C_c^\infty$ function of  $(x, \xi)$, and for simplicity we take it to be smooth in $h \in [0, h_0]$. 
By abuse of notation we denote the principal symbol $p(x, \xi, 0)$ by $p(x, \xi)$, and we will write $p(x, hD)$ for $p(x, hD, h)$. 

Following \cite{koch07}, we assume that our family of quasimodes $p(x,hD)u(h)=O_{L^{2}}(h)$ is semiclassically localised:

\begin{defin}\label{localised}
A function $u$ depending parametrically on $h$ is localised if there exists $\chi\in{}C_{c}^{\infty}(T^{\star}M)$ such that
$$u=\chi(x,hD)u+O_{\mathcal{S}}(h^{\infty})$$
where $\mathcal{S}$ is the space of Schwartz functions, and $ g \in O_{\mathcal{S}}(h^{\infty})$ means that each seminorm of $g$ is $O(h^\infty)$.
\end{defin}

Localisation is compatible with the assumption that $p(x,hD)u=O_{L^{2}}(h)$: that is, if $\chi\in{}C_{c}^{\infty}(T^{\star}M)$ then
$$p(x,hD)u=O_{L^{2}}(h)\Rightarrow{}p(x,hD)(\chi(x,hD)u)=O_{L^{2}}(h).$$
Using this localisation assumption we are able to turn the global problem into a local problem on small patches in $T^{\star}M$. If $\chi\in{}C_{c}(T^{\star}M)$ such that
$$u=\chi(x,hD)u+O_{\mathcal{S}}(h^{\infty})$$ then, using compactness of the support of $\chi$, we can write
$$\chi(x,\xi)=\sum_{i=1}^{N}\chi_{i}(x,\xi)$$
for some $N<\infty$ where each $\chi_{i}$ has arbitrarily small support. In this fashion we reduce estimating $\norm{\chi(x,hD)u}_{L^{p}(H)}$ to estimates on $\norm{\chi_{i}(x,hD)u}_{L^{p}(H)}$ (the error term $O_{\mathcal{S}}(h^{\infty})$ is of course trivial to estimate). Due to this localisation we can replace $M$ with $\R^{n}$, $H$ with $\R^{n-1}$ and $T^{\star}M$ with $\R^{n}\times{}\R^{n}$. We write $x\in{}M$ as $x=(y,\bdf)$ where $y\in{}H$ and $\bdf$ is the normal direction to $H$.

Still following \cite{koch07}, we  further reduce this problem to localising around points $(x_{0},\xi_{0})$ where $p(x_{0},\xi_{0})=0$. To achieve this we use Lemma 2.1 of \cite{koch07} which shows that if $|p(x,\xi)|\geq{1/C}$ on a local patch then we can invert $p(x,hD)$ up to order $h^{\infty}$. That is, choosing $\chi(x,\xi)$ supported on this patch, we can find some $q(x,hD)$ such that
 $$q(x,hD)p(x,hD)\chi(x,hD)=\chi(x,hD)+O_{L^{2}\rightarrow{}L^{2}}(h^{\infty})$$
 and
 $$p(x,hD)q(x,hD)\chi(x,hD)=\chi(x,hD)+O_{L^{2}\rightarrow{}L^{2}}(h^{\infty}).$$
 So if $p(x,hD)u=O_{L^{2}}(h)$ and $|p(x,\xi)|>1/C$ we can invert $p(x,hD)$ to get
$$\chi(x,hD)u=O_{L^{2}}(h).$$

We can combine this estimate with the following `semiclassical Sobolev inequality' (see \cite{burq2}, \cite{evans} or \cite{koch07} for proof) to obtain hypersurface restriction estimates.

\begin{lem}[semiclassical Sobolev estimates]\label{semiclassicalsobolev}
Suppose that a family $u=u(h)$ satisfies the localisation condition. Then for $1\leq{}q\leq{}p\leq\infty$
$$\norm{u}_{L^{p}}\lesssim{}h^{n(1/p-1/q)}\norm{u}_{L^{q}}+O(h^{\infty}).$$

\end{lem}

To get the $L^{2}$ norm of the restriction of $u$ to $H$ we use Lemma \ref{semiclassicalsobolev}  in  only the $\bdf$ coordinates. This is justified as localisation in $T^{\star}\R^{n}$ implies localisation in $T^{\star}\R^{n-1}$ (see \cite{tacy09}). We have
\begin{equation}\norm{u(y,0)}_{L^{2}_{y}}\lesssim\norm{u(y,\bdf)}_{L^{\infty}_{\bdf}L^{2}_{y}}\lesssim{}h^{-\frac{1}{2}}\norm{u(y,z)}_{L^{2}_{z}L^{2}_{y}}.\label{goodL2}\end{equation}
So, if $|p(x,\xi)|\geq{1/C}$, and $Pu = O(h)$, the $L^{2}$ norm of $u$ when restricted to a hypersurface $H$ containing $x_0$ is $O(h^{\frac{1}{2}})$. This is significantly better than the $L^{2}$ estimate given by Theorem \ref{glancingflowtheorem}. Consequently we can ignore regions where $p(x,\xi)$ is bounded away from zero.

To get better estimates when $p(x_0, \xi_0) = 0$ than what can be obtained from Lemma~\ref{semiclassicalsobolev} (which uses only localisation), we need to make assumptions on the function $p$ (to prevent $p$ vanishing identically, for example, in which case the assumption $Pu = O(h)$ is vacuous!). Our first assumption (A1) is that $p$ vanishes simply on each cotangent fibre:

 \begin{itemize}
 \item[(A1)] 
 \text{ for any point $(x_{0},\xi_{0})$ such that $p(x_{0},\xi_{0})=0$, $\partial_{\xi}p(x_{0},\xi_{0})\neq{}0$.}
 \end{itemize}
 
 Our second condition is a geometric condition on the characteristic variety. The condition eliminates examples such as $p(x, \xi) = \xi_1$, i.e. $P = h D_{x_1}$, for which we cannot estimate $\| u \|_{L^2(H)}$ by better than the $h^{-1/2}$ estimate given by Lemma~\ref{semiclassicalsobolev} alone. Let us note that (A1) implies that the set 
\begin{equation}
\{\xi\mid{}p(x_{0},\xi)=0\}\subset{}T_{x_{0}}^{\star}M
\label{charvar}\end{equation}
 is a smooth hypersurface in $T_{x_{0}}^{\star}M$.

\begin{itemize}
\item[(A2)] 
\text{For each $x_0 \in M$, the second fundamental form of \eqref{charvar} is positive definite.}
\end{itemize}

\begin{defin}\label{admissible}
A symbol $p(x,\xi)$ is \emph{admissible} if it satisfies condition (A1) and (A2). \end{defin}

In addition we make the geometric assumption of curvature with respect to the flow.

\begin{defin}\label{curvedtoflow}
A hypersurface $H$ of $M$ is  \emph{curved} with respect to the flow if the projection of the bicharacteristic flow to $M$ is at most simply tangent to $H$, or in other words, if for one (and hence any) boundary defining function $\bdf$ for $H$, we have
\end{defin}
\begin{itemize}
\item[(A3)] \text{ For any $(x_0, \xi_0)$, $\dot \bdf(x_0, \xi_0) = 0$ implies that $\ddot \bdf (x_0, \xi_0) \neq 0$.}
\end{itemize} 

\begin{remark} In the case $P(h) = h^2 \Delta - 1$, where $\Delta$ is the Laplacian on $M$ with respect to a Riemannian metric, assumptions (A1) and (A2) are satisfied, and (A3) is satisfied iff $H$ has positive definite second fundamental form. Thus, in this case our result reduces to that of Burq-G\'erard-Tzvetkov \cite{burq07} ($n=2$) and Hu \cite{hu09} ($n \geq 2$). 
\end{remark}


\section{Evolution equation}\label{symbolfac}

Using the argument in the previous section we can assume that $p(x_0, \xi_0) = 0$. 
Assumption (A1) then tells us that  $\partial_{\xi}p(x_{0},\xi_{0})\neq{}0$. Let us choose coordinates $x = (y, \bdf)$ where $y \in \RR^{n-1}$ and $\bdf \in \RR$ is a boundary defining function for $H$. Let $\xi = (\eta, \nu)$ be the dual coordinates. If $\partial_\nu p(x_0, \xi_0) \neq 0$ then we have $\dot r \neq 0$ and, as mentioned in the Introduction
(see \eqref{nonconcentration}), $u$ does not concentrate at $H$ at all. 
Therefore we may assume that $\partial_\nu p(x_0, \xi_0) = 0$.
Therefore we have $\partial_\eta p(x_0, \xi_0) \neq 0$. By a linear change of $y$ coordinates we can assume that $\partial_{\eta_1} p(x_0, \xi_0) \neq 0$ and $\partial_{\eta_j} p(x_0, \xi_0) = 0$ for $j \geq 2$. 

Now we apply the implicit function theorem and deduce that the characteristic variety $\{ p = 0 \}$ implicitly defines $\xi_1$ as a  smooth function of $(x, \xi_2, \dots, \xi_n)$:
\begin{equation}
p = 0 \implies \xi_1 = a(x, \xi_2, \dots, \xi_n).
\label{adefn}\end{equation}
We shall now write $x_1 = t$ and think of it as a time variable. We write $x = (t, \xbar)$ and similarly, $\xi_1 = \tau$ and $\xi = (\tau, \xibar)$. We also write $y = (t, y')$ and $\eta = (\tau, \eta')$. Thus $x = (t, y', r)$ and correspondingly $\xi = (\tau, \eta', \nu)$. 
We write $T$ for the `initial' hypersurface $ \{ t = 0 \}$, and recall that $H = \{ r = 0 \}$. We assume that $t = 0$ at $(x_0, \xi_0)$ and write
$(x_0, \xi_0) = ((0, \xbar_0), \xi_0) = ((0, y_0', 0), (\tau_0, \eta'_0, \nu_0))$. 

As a consequence of \eqref{adefn}, we have 
$$p = e(x, \xi) \big(\tau - a(x, \xibar)\big)$$
near $(x_0, \xi_0)$, where $e(x_0, \xi_0) \neq 0$. By localising suitably we may assume that $e \neq 0$ on the support of our localising function $\chi$. The condition $Pu = O(h)$ in $L^2$ then implies that
$$
e(x, hD_x) \big(h D_t - a(x, hD_{\xbar}) \big) u = O_{L^2(M)}(h)
$$
and using the local invertibility modulo $O(h^\infty)$ of $e(x, hD_x)$, we find that 
\begin{equation}\big(hD_{t}-a(x,hD_{\xbar}) \big)u=hf(t,\xbar)
\label{evolutionequation}\end{equation}
where $\norm{f}_{L^{2}(M)}=O_{L^{2}}(1)$.

We view \eqref{evolutionequation} as an evolution equation for $u$, which determines $u$ given the `initial data' $u(0, \xbar)$ and the
inhomogeneous term $f(t, \xbar)$. This determines a family of solution operators $U_{s}(t)$, such that $U_{s}(t)$ is the solution operator for the evolution equation
$$
\big(hD_{t}-a(s+t,\xbar,hD_{\xbar}) \big)u=0, \quad u(0, \xbar ) = u(\xbar)
$$

Using Duhamel's principle we write

$$u(t,\bar{x})=U_{0}(t)u(0,\bar{x})+i\int_{0}^{t}U_{s}(t-s)f(s,\bar{x})ds.$$
Now let $R_H$ be the operation of restriction to the hypersurface $H$, and let $W_{s}(t) = R_H \circ U_{s}(t)$. Also, let $u_{0}=u(0,\xbar)$ be the restriction of $u$ to the initial hypersurface $T = \{ t = 0 \}$.
We then have
$$u(t,y',0)= W_{0}(t)u_0 +i\int_{0}^{t}W_{s}(t-s)f(s,\xbar)ds.$$
Using Minkowski's inequality we have 
\begin{multline}
\norm{u}_{L^{2}(H)}\lesssim\left(\int\norm{W_{0}(t)u_{0}}_{L^{2}_{y'}}^{2}dt\right)^{1/2}+\\
\int_{\R}\left(\int\norm{W_{s}(t-s)f(s,\xbar)}_{L^{2}_{y'}}^{2}dt\right)^{1/2}ds
\label{qnorm}\end{multline}
We recall from \eqref{nonconcentration} (with $H' = T$) that $\| u_0 \|_{L^2(T)} \lesssim \| u \|_{L^2(M)}$. 
Therefore, to prove Theorem~\ref{glancingflowtheorem}, i.e.  obtain a $L^{2}$ bound of
$$\norm{u}_{L^{2}(H)}\lesssim{}h^{-1/6}\norm{u}_{L^{2}(M)}$$
it suffices to obtain an estimate, uniform in $s$, of the form
\begin{equation}\left(\int\norm{W_{s}(t-s) f}_{L_{y'}^{2}}^{2}dt\right)^{1/2}\lesssim{}h^{-1/6}\norm{f}_{L^{2}(T)}.\label{Wbounds}\end{equation}
For each $s$ we will show that \eqref{Wbounds} holds with a constant that depends only on the seminorms of $a(x,\xibar)$. 
In fact, the estimates are uniform given uniform bounds on a finite number of derivatives of $a$, and given uniform lower bounds on the nondegeneracies involved in the computation in Section~\ref{canonical} --- see Remark~\ref{comech}. Such uniform bounds hold provided that the patch size is chosen sufficiently small. Therefore we only address the estimate for $W_{0}(t)$, which we denote by $W(t)$ from here on. 
To obtain this estimate we view $W(t)$, thought of as a single operator from $L^2(T)$ to $L^2(H)$ instead of as a family parametrised by $t$, as an Fourier integral operator. 


\section{Fourier integral representation}\label{oscillatoryevol}
We need to express the solution operator for the evolution equation
\begin{equation}hD_{t}-a(t,\bar{x},hD_{\bar{x}})=0\label{evolution}\end{equation}
as an Fourier integral operator.  We will then be able to transfer properties of the flow to properties of the phase function defining the operator $U(t)$. 

\begin{prop}\label{FIO}
Suppose $U(t):L^{2}(\R^{d})\rightarrow{}L^{2}(\R^{d})$ satisfies 
$$hD_{t}U(t)-A(t)U(t)=0,\quad{}U(0)=Id$$
where A(t) is a semiclassical pseudodifferential operator such that the symbol $a(t,\bar{x},\eta)$ of $A(t)$ is  real and is smooth in $h$. Then there exists some $t_{0}>0$ independent of $h$ such that for $0\leq{}t\leq{}t_{0}$
$$U(t)u(\bar{x})=\frac{1}{(2\pi{}h)^{d}}\int\int{}e^{\frac{i}{h}(\phi(t,\bar{x},\xibar)-\wbar\cdot\xibar)}b(t,\bar{x},\xibar,h)u(\wbar)d\wbar d\xibar+E(t)u(\bar{x})$$
where
$$\partial_{t}\phi(t,\bar{x},\xibar)-a(t,\bar{x},\partial_{\bar{x}}\phi(t,\bar{x},\xibar))=0,\quad{}\phi(0,\bar{x},\xibar)=\bar{x}\cdot\xibar$$
$$b(t,\bar{x},\xibar,h)\in{}C^{\infty}_{c}(\R\times{}T^{\star}\R^{d}\times{}\R)\quad{}E(t)=O(h^{\infty}):S'\rightarrow{}S$$

\end{prop}
\begin{proof}
This is in fact the normal parametrix construction yielding the eikonal equation for the phase function. See \cite{evans} Section 10.2.
\end{proof}

Recall that $W(t)=R_{H}\circ{}U(t)$ so we have 
$$
W(t)f(y')=\frac{1}{(2\pi{}h)^{n-1}}\iint{}e^{\frac{i}{h}(\phi(t,(y',0),\xibar)-\wbar\cdot\xibar)}b(t,y',\eta,h)f(\wbar)d\wbar d\xibar
$$
In what follows we will write $\phi(t,y', \eta', \nu)$ for $\phi(t,(y', 0),\xibar)$ (recall that $\xibar = (\eta', \nu)$). We want to estimate the operator norm of $W(t)$ regarded as a single operator acting from $L^2(T)$ to $L^2(H)$. Note that $W(t)=Z \circ \mathcal{F}_h$ where $\mathcal{F}_h$ is the semiclassical Fourier transform:
$$\mathcal{F}_h f(\xibar)=\frac{1}{(2\pi{}h)^{\frac{n-1}{2}}}\int{}e^{-\frac{i}{h}\xibar\cdot{}\vbar}f(\vbar)d\vbar$$
and the operator $Z$ is given by 
$$
Z g(t, y') = \frac{1}{(2\pi{}h)^{\frac{n-1}{2}}}\iint{}e^{\frac{i}{h}\phi(t, y', \eta', \nu)}b(t,y',\eta', \nu,h) g(\eta', \nu) \, d\eta' \, d\nu.
$$
As $\norm{\mathcal{F}_{h}f}_{L^{2}}=\norm{f}_{L^{2}}$ it is enough to estimates $L^{2}\to{}L^{2}$ operator norm of $Z$. To estimate the operator norm of $Z$ we view it as a semiclassical Fourier integral operator and analyse its canonical relation. 


\section{Canonical relation}\label{canonical}
To prove Theorem~\ref{glancingflowtheorem} we need to show that the operator norm of $Z$ is bounded by $Ch^{-1/6}$. To do this we use the following theorem of Pan and Sogge \cite{pan90} which is the analogue for oscillatory integral operators of Melrose and Talyor's \cite{melrose85} theorem on Fourier integral operators with folding canonical relations. 

\begin{thm}\label{foldingrelations} 
Let the oscillatory integral operator $T_{\lambda}$ be defined by
$$T_{\lambda}f(x)=\int_{\R^{d}}e^{i\lambda\psi(x,y)}\beta(x,y)f(y)dy$$
where $\beta\in{}C_{0}^{\infty}(\R^{d}\times{}\R^{d})$ and the phase function $\psi\in{}C^{\infty}(\R^{d}\times{}\R^{d})$ is real. If the left and right projections from the associated canonical relation
$$\canon_{\psi}=\{(x,\psi'_{x}(x,y),y,-\psi'_{y}(x,y))\}$$ 
are at most folding singularities then
$$\norm{T_{\lambda}f}_{L^{2}(\R^{d})}\lesssim{}\lambda^{-\frac{d}{2}+1/6}\norm{f}_{L^{2}(\R^{d})}$$
\end{thm}

Let us recall (see for example \cite{guillemingolubitsky}) that a smooth map $F : \RR^d \to \RR^d$ has a folding singularity
at $x \in \RR^d$ if
\begin{itemize}
\item[(i)] $dF(x)$ is rank $d-1$, 
\item[(ii)] the function $\det dF$ vanishes simply at $x$, implying in particular that locally near $x$, the set of $y \in \RR^d$ such that $dF(y)$ has rank $d-1$ is a smooth hypersurface $S$ containing $x$, and \item[(iii)] the kernel of $dF(x)$ is not contained in the tangent space to $S$:
$$
T_x S + \ker dF(x) = T_x \RR^d.
$$
\end{itemize}
Given (i) an equivalent condition to (ii) and (iii) is that, if $v$ is a nonzero element of $\ker dF(x)$, then 
\begin{equation}
D_v(\det dF(x)) \neq 0.
\label{iii}\end{equation}

The operator $Z$ is a Fourier integral operator with canonical relation
$$
\canon = \big\{ (t, y', \partial_t \phi, \partial_{y'} \phi, \eta', \nu, -\partial_{\eta'} \phi, -\partial_{\nu} \phi ) \big\}.
$$
The left and right projections on $\canon$ are represented in local coordinates by
$$
\pi_L : (t, y', \eta', \nu) \mapsto (t, y', \partial_t \phi, \partial_{y'} \phi)
$$
and 
$$
\pi_R :  (t, y', \eta', \nu) \mapsto (\eta', \nu, \partial_{\eta'} \phi, \partial_{\nu} \phi )
$$
(where we removed the irrelevant minus signs from $\pi_R$ for notational convenience). 

 The matrix $d\pi_L$ takes the form

$$d\pi_{L}=\left(\begin{BMAT}(e)[2pt,1.5cm,1.5cm]{c.c}{c.c}
 \Id & 0 \\
 * & B \end{BMAT}\right)
$$
where
$$
B=\left(\begin{BMAT}(e){ccc.c}{c.ct}
&\partial^{2}_{t\eta'}\phi& & \partial^{2}_{t\nu}\phi \\ 
&&&\\
 &\partial^{2}_{y'\eta'}\phi& & \partial^{2}_{y'\nu}\phi \end{BMAT}\right)
$$
At $(x_0, \xi_0)$ we have $\partial^2_{y' \eta'} \phi = \Id$, $\partial^2_{t\eta'} \phi = \partial_{\eta'} a = 0$, $\partial^2_{y' \nu} \phi = 0$ and $\partial^2_{t \nu} \phi = \partial_{\nu} a = 0$, so we get 
$$
B =\left(\begin{BMAT}(e){ccc.c}{c.ct}
&0& & 0 \\ 
&&&\\
 &\Id& & 0 \end{BMAT}\right).
$$
It is clear that the vector field $\partial_\nu$ is in the kernel of $d\pi_L(x_0, \xi_0)$. Moreover, $\det d\pi_L$ is given by $\partial^2_{t\nu} \phi \cdot \det (\partial^2_{y' \eta'} \phi)$ plus terms vanishing to second order at $(x_0, \xi_0)$. To show that $\pi_L$ has a fold at $(x_0, \xi_0)$ we need by \eqref{iii} to show that $\partial_\nu (\det d\pi_L)$ is nonzero at $(x_0, \xi_0)$. 
Due to the vanishing of both `off-diagonal' terms $\partial^2_{t\eta'} \phi$ and $\partial^2_{y' \nu} \phi$, the nonvanishing of $\partial_\nu (\det d\pi_L)$ at 
$(x_0, \xi_0)$ is equivalent to the nonvanishing of $\partial_{\nu} (\partial^2_{t\nu} \phi) = \partial^3_{t\nu \nu} \phi$.

The matrix $d\pi_R$ takes the form 
$$
d\pi_R=\left(\begin{BMAT}(e)[2pt,1.5cm,1.5cm]{c.c}{c.c}
              0 & \Id \\
	D	& *
             \end{BMAT}\right)$$

$$
D=\left(\begin{BMAT}(e){c.ccc}{bc.c}
 \partial^{2}_{\eta't}\phi& & \partial^{2}_{y'\eta'}\phi &\\
&&&\\
\partial^{2}_{\nu{}t}\phi& &\partial^{2}_{\nu{}y'}\phi &\end{BMAT}\right)
$$
and we see that $\partial_t$ is in the kernel of $d\pi_R(x_0, \xi_0)$. 
To show that $\pi_R$ has a fold at $(x_0, \xi_0)$ we need by \eqref{iii} to show that $\partial_t (\det d\pi_L)$ is nonzero at $(x_0, \xi_0)$. 
As above, due to the vanishing of the `off-diagonal' terms $\partial^2_{t\eta'} \phi$ and $\partial^2_{y' \nu} \phi$, the nonvanishing of $\partial_t (\det d\pi_L)$ at 
$(x_0, \xi_0)$ is equivalent to the nonvanishing of $\partial_t (\partial^2_{t\nu} \phi) = \partial^3_{t t \nu} \phi$. 

The proof of Theorem~\ref{glancingflowtheorem} is therefore completed by the following Lemma:

\begin{lem} Under assumptions (A1), (A2) and (A3), we have
$$
\partial^3_{t \nu \nu} \phi(x_0, \xi_0) \neq 0, \text{ and } \partial^3_{t t \nu} \phi(x_0, \xi_0) \neq 0.
$$
\end{lem}

\begin{remark} To simplify notation we write $(x_0, \xi_0)$ for the argument of $\phi$ corresponding to this point, although $(0, y_0',0,\tau_{0}, \eta_0', \nu_{0})$ would be more accurate. 
\end{remark}

\begin{proof} 
We use the Hamilton-Jacobi equation 
\begin{equation}
\partial_{t}\phi(t, \xbar, \eta', \nu) = a\big(t, \xbar,\partial_{\xbar}\phi(t,  \xbar,\eta', \nu)\big).
\label{HJ}\end{equation}
Since at $t=0$ we have $\phi(0, \xbar, \eta', \nu) = y' \cdot \eta' + \bdf \nu$ (recall that $\xbar = (y', \bdf)$), we have 
$$
\partial^3_{t \nu \nu} \phi(0,\xbar,\eta', \nu) = \partial^2_{\nu \nu} a(0, \xbar, \eta', \nu).
$$
Now we apply assumption (A2): it says that the second fundamental form of the submanifold $\{ \tau = a(x, \xibar_0) \} \subset T_{x_0} M$ is positive definite. Since $\partial_{\xibar} a = 0$ at $(x_0, \xibar_0)$, the second fundamental form of this submanifold at $(x_0, \xi_0)$ is given by the matrix of second derivatives of $a$:
$$
h_{ij}(\xi_0) = \partial^2_{\xi_i \xi_j} a(x_0, \xibar_0), \quad 2 \leq i, j \leq n.
$$
Therefore, $\partial^2_{\nu \nu} a \neq 0$ at $(x_0, \xibar_0)$, showing that $\pi_L$ has a fold singularity at $(x_0, \xi_0)$.

To treat the term $\partial^3_{t t \nu} \phi(x_0, \xi_0)$, we differentiate \eqref{HJ} in $t$, obtaining
$$
\partial^2_{tt} \phi = \partial_t a + \partial_{\xibar} a\cdot \partial^2_{\xbar t} \phi.
$$
Using \eqref{HJ} again on the term $\partial^2_{\xbar t} \phi$ we obtain
$$
\partial^2_{tt} \phi = \partial_t a + \partial_{\xibar} a \cdot\big( \partial_{\xbar} a + \partial_{\xibar} a \cdot\partial^2_{\xbar \xbar} \phi \big).
$$
We evaluate this at $t=0$ since the next derivative to be applied, namely $\partial_\nu$, is tangent to $\{ t = 0 \}$. At $t=0$, we have $\partial^2_{\xbar \xbar} \phi = 0$, so we get 
$$
\partial^2_{tt} \phi  \Big|_{t=0} = \partial_t a + \partial_{\xibar} a \cdot \partial_{\xbar} a.
$$
Now when we differentiate in $\nu$, we get
$$
\partial^3_{tt\nu} \phi(x_0, \xi_0) = \partial^2_{t \nu} a(x_0, \xibar_0) + \partial^2_{\xibar \nu} a(x_0, \xibar_0)\cdot \partial_{\xbar} a(x_0, \xibar_0)
$$
since $\partial_{\xibar} a(x_0, \xibar_0) = 0$. 

At this point we remind the reader that we have chosen coordinate $(t,y',\bdf)$ and $(\tau,\eta',\nu)=(\tau,\xibar)$ such that
$$\partial_{\xibar}p(x_{0},\xi_{0})=0$$
and
$$\tau_{0}-a(t_{0},\xbar_{0},\xibar_{0})=0.$$
It follows that
\begin{equation}\begin{cases}
  \partial_{\xi}a(x_{0},\xibar_{0})=0\\
\partial_{\tau}p(x_{0}.\xi_{0})=e(x_{0},\xi_{0})\\ 
\partial_{\bar{x}}p(x_{0},\xi_{0})=-e(x_{0},\xi_{0})\partial_{\bar{x}}a(x_{0},\xibar_{0})\\
\partial_{t}p(x_{0},\xi_{0})=-e(x_{0},\xi_{0})\partial_{t}a(x_{0},\xi_{0}).
  \end{cases}
\label{miscid}\end{equation}

Now we apply assumption (A3), which says that $\ddot \bdf \neq 0$. We express $\ddot \bdf$ in terms of $a$. We have 
$$
\dot \bdf = \partial_\nu p = \partial_\nu \big(e (\tau - a) \big).
$$
Differentiating a second time and using the flow identities
$$\dot{x}=\partial_{\xi}p(x,\xi)\quad\dot{\xi}=-\partial_{x}p(x,\xi),$$
we have
$$
\ddot \bdf = \Big( \partial_\tau p \partial_t + \partial_{\xibar} p \partial_{\xbar} - \partial_t p \partial_\tau - \partial_{\xbar} p \partial_{\xibar} \Big) \Big(  (\tau - a)\partial_\nu e - e \partial_\nu a \Big).
$$
At $(x_{0},\xi_{0})$ using the identities given in \eqref{miscid} we can simplify this to
$$
\ddot \bdf(x_0, \xi_0) = -e \Big( \partial^2_{\nu t} a(x_0, \xibar_0) + \partial_{\xbar} a(x_0, \xibar_0)\cdot \partial^2_{\nu \xibar} a(x_0, \xibar_0) \Big).
$$
Therefore, applying assumption (A3), we find
$$
\partial^3_{tt\nu} \phi(x_0, \xi_0) = - \frac{\ddot \bdf(x_0, \xi_0)}{e(x_0, \xi_0)} \neq 0.
$$
This shows that $\pi_R$ has a fold singularity at $(x_0, \xi_0)$ and completes the proof. 
\end{proof}

\begin{remark} It is easy to see from the calculations above that 
assumption (A3) is \emph{equivalent} to the statement that $\pi_R$ has a folding singularity. Similarly,  
assumption (A2) is equivalent to the statement that $\pi_L$ has a folding singularity for every hypersurface $H$ whose tangent space $T_{x_0} H$ at $x_0$ contains $\partial_\xi p(x_0, \xi_0) \partial_{x}$, 
i.e. the tangent vector of the projected bicharacteristic through $(x_0, \xi_0)$.  
\end{remark}

\begin{remark}\label{comech}
According to \cite{comech}, Theorem 2.2, one obtains uniform bounds of the form $C h^{-1/6}$ on the norms of the operators $W_s$ given by \eqref{Wbounds} provided that there are uniform bounds on a finite number of derivatives of the symbol of $W_s$, and uniform lower bounds on the determinant of $\partial^2_{y' \eta'} \phi$, $\partial_\nu (\partial^2_{t \nu} \phi)$, and $\partial_t (\partial^2_{t \nu} \phi)$. These lower bounds are achieved simply by shrinking the patch size sufficiently and using continuity. Thus we obtain a bound as in \eqref{Wbounds} uniformly in $s$, as desired. 
\end{remark}


\section{Optimality of Theorem~\ref{glancingflowtheorem}}

All the estimates given by Thereom \ref{glancingflowtheorem} are sharp. We study a simple local model around $(0,0)$ for hypersurface curved with respect to the flow. Let $H=\{x\mid{}x_{n}=0\}$ and $p(x,\xi)$ be given by
$$p(x,\xi)=\xi_{1}-x_{n}-\sum_{i=2}^{n}\xi_{i}^{2}$$
Note that
\begin{align*} \dot{t}&=1&\dot{y'}&=-2\eta'&\dot{r}&=-2{}\nu
\\
\dot{\tau}&=0&\dot{\eta'}&=0&\dot{\nu}&=1
\end{align*}
Therefore the flow $(x(s),\xi(s))$ with intial point $(0,0)$ is given by 
\begin{align*}
t(s)&=s&{}y'(s)&=0&r(s)&=-s^{2}\\
\tau(s)&=0&\eta'(s)&=0&\nu(s)&=s
\end{align*}
So we have that condition (A3) is clearly satified as $\ddot{r}(0)=-2$. 
$$p(x,hD)=hD_{t}-r-h^{2}D_{r}^{2}-\sum_{i=1}^{n-2}h^{2}D_{y'}^{2}$$
 It is easier to develop a solution in Fourier space. Note that
$$
\mathcal{F}_{h} \circ p(x,hD) \circ \mathcal{F}_{h}^{-1} =\tau-hD_{\nu}-\nu^{2}-\eta'\cdot\eta'.
$$
As the semiclassical Fourier transform preserves $L^{2}$ norms if
$$\norm{(\tau-hD_{\nu}-\nu^{2}-\eta'\cdot\eta')f}_{L^{2}}=O_{L^{2}}(h)$$
and $u=\mathcal{F}_{h}^{-1}f$, then
$$\norm{p(x,hD)u}_{L^{2}}=O_{L^{2}}(h).$$ 
We therefore seek a solution for
\begin{equation}(\tau-hD_{\nu} - \nu^2 - \eta'\cdot\eta')f=0;\label{Fdiffeq}\end{equation}
it is obvious that
$$g(\tau,\eta',\nu)=e^{\frac{i}{h}(\frac{1}{3}\nu^{3}+\nu(\tau-\eta'\cdot\eta'))}$$
is a solution to \eqref{Fdiffeq}. The natural scaling $\nu\to{}h^{-1/3}\nu$ induces a scaling of $\tau\to{}h^{-2/3}\tau$ and $\eta'\to{}h^{-1/3}\eta'$ and  accordingly we place cut off functions appropriate to that scale. Let
$$f(\tau,\eta',\nu)=h^{-\frac{n-2}{6}-\frac{1}{3}}\chi(|\nu|)\chi(h^{-\frac{2}{3}}|\tau|)\chi(h^{-\frac{1}{3}}|\eta'|)e^{\frac{i}{h}\psi(\tau,\eta',\nu)}$$
where
$$\psi(\tau,\eta',\nu)=\frac{1}{3}\nu^{3}+\nu(\tau-\eta'\cdot\eta')$$
Now $\norm{f}_{L^{2}}=O_{L^{2}}(1)$ and  $f$ satisfies \eqref{Fdiffeq} up to an $O(h)$ error coming from the $D_\nu$ hitting the cutoff function $\chi(|\nu|)$.  We define the function $u$ as
$$u=\chi(|x|)\mathcal{F}_{h}^{-1}f$$
Now $R_{H}u$ is given by
$$R_{H}u(y)=\frac{h^{-\frac{n-2}{6}-\frac{1}{3}}\chi(|y|)}{(2\pi{}h)^{\frac{n}{2}}}\int{}e^{\frac{i}{h}(t\tau+y' \cdot \eta'+\psi(\tau,\eta',\nu))}\chi(|\nu|)\chi \left( \frac{|\tau|}{h^{2/3}} \right) \chi \left( \frac{|\eta'|}{h^{1/3}} \right)d\tau{}d\nu{}d\eta' .$$
For $|t|\leq{} \epsilon h^{1/3}$, $\epsilon$ small, the factor $e^{\frac{i}{h}t\tau}$ does not oscillate significantly and can be ignored. Similarly for $|y'|\leq{} \epsilon h^{2/3}$ the factor $e^{\frac{i}{h}y' \cdot \eta'}$ does not oscillate significantly and is also ignored. On the other hand, there are oscillations in the $\nu$ variable. At $\tau = \eta' = 0$ there is degenerate stationary phase at $\nu = 0$; Theorem 7.7.18 of \cite{Hormandervol1} applies and shows that there is a lower bound of the form
$$|R_{H}u(t, y')| \sim {}h^{-\frac{n-2}{6}-\frac{1}{3}-\frac{n}{2}+\frac{n-1}{3}+\frac{2}{3}}=h^{-\frac{n-1}{3}}$$
for $|t|\leq{} \epsilon h^{1/3}$,  $|y'|\leq{} \epsilon h^{2/3}$. Thus on this set we get a lower bound on the $L^p$ norm:
$$\norm{u}_{L^{p}([0,\epsilon h^{1/3}]_t\times{}B(0,\epsilon h^{2/3})_{y'})}\sim{}h^{-\frac{n-1}{3}+\frac{1}{3p}+\frac{2(n-2)}{3p}}=h^{-\left(\frac{n-1}{3}-\frac{2n-3}{3p}\right)}$$
which saturates the estimate of Theorem \ref{glancingflowtheorem}.

\bibliography{glancingflow}
\bibliographystyle{plain}

\end{document}